\documentclass{gtart_h}  

\def\ifplaintex{\expandafter\ifx\csname documentclass\endcsname\relax}

\ifplaintex 
\else
\fi

\expandafter\ifx\csname beginpicture\endcsname\relax
\expandafter\ifx\csname documentclass\endcsname\relax
\input pictex \else
\input prepictex \input pictex \input postpictex \fi\fi

\def\gt{{\mathsurround=0pt\it $\cal G\mskip-2mu$eometry \&\ 
$\cal T\!\!$opology}}        

\def\gtp{{\mathsurround=0pt\it $\cal G\mskip-2mu$eometry \&\ 
$\cal T\!\!$opology $\cal P\!$ublications}}  

\def\lognumber#1{\def\thelognumber{#1}}
\def\volumenumber#1{\def\thevolumenumber{#1}}
\def\papernumber#1{\def\thepapernumber{#1}}
\def\volumeyear#1{\def\thevolumeyear{#1}}

\def\pagenumbers#1#2{\def\startpage{#1}\def\finishpage{#2}}
\def\published#1{\def\publishdate{#1}}
\def\proposed#1{\def\theproposer{#1}}
\def\seconded#1{\def\theseconders{#1}}
\def\received#1{\def\receiveddate{#1}}
\def\revised#1{\def\reviseddate{#1}}
\def\accepted#1{\def\accepteddate{#1}}
\def\asciititle#1{\def\theasciititle{#1}}

\long\def\asciiabstract#1{\long\def\theasciiabstract{#1}}
\def\asciikeywords#1{\def\theasciikeywords{#1}}

\let\\\par\let\thelognumber\relax
\let\thevolumenumber\relax\let\thepapernumber\relax
\let\thevolumeyear\relax\let\thesamplenumber\relax\let\startpage\relax
\let\finishpage\relax\let\publishdate\relax\let\receiveddate\relax
\let\reviseddate\relax\let\accepteddate\relax\let\theasciititle\relax
\let\theasciiauthors\relax
\let\theasciiabstract\relax\let\theasciikeywords\relax
\let\theasciiemail\relax\let\theshortauthors\relax\let\theshorttitle\relax

\long\def\maketitlep{   

\count0=\startpage

\gt\hfill      
\beginpicture
\setcoordinatesystem units <0.33truein, 0.33truein> point at 2.2 0.9
\setplotsymbol ({$\cal G$})
\plotsymbolspacing=9truept
\circulararc 315 degrees from 0 1 center at 0 0
\setplotsymbol ({$\cal T$})
\circulararc 315 degrees from 1 -1 center at 1 0
\endpicture
\break
{\small\ifx\thesamplenumber\relax 
Volume \else Sample
\fi\thevolumenumber\ (\thevolumeyear)
\startpage--\finishpage\nl
Published: \publishdate}
\vglue 0.5truein plus 0.4fil minus 0.1truein

{\parskip=0pt\leftskip 0pt plus 1fil\def\\{\par\smallskip}{\ifplaintex\large
\else\Large\fi\bf\thetitle}\par\medskip}   

\vglue 0pt plus 0.1fil 

{\parskip=0pt\leftskip 0pt plus 1fil\def\\{\par}{\sc\theauthors}
\par\medskip}

\vglue 0pt plus 0.1fil 

{\small\parskip=0pt\let\newline\\
{\leftskip 0pt plus 1fil\def\\{\par}{\sl\theaddress}\par}
\expandafter\ifx\theemail\relax    
\relax\else\vglue 5pt plus 0.02fil minus 2pt\def\\{\stdspace{\rm 
and}\stdspace} 
\cl{Email:\stdspace\tt\theemail}\fi
\ifx\theurl\relax                  
\relax\else\vglue 5pt plus 0.02fil minus 2pt\def\\{\stdspace{\rm 
and}\stdspace}
\cl{URL:\stdspace\tt\theurl}\fi\par}

\vglue 7pt plus 0.3fil minus 3pt

{\bf Abstract}
\vglue 5pt plus 0.1fil minus 2pt

\theabstract

\vglue 7pt plus 0.3fil minus 3pt

{\bf AMS Classification numbers}\quad Primary:\quad \theprimaryclass

Secondary:\quad \thesecondaryclass

\vglue 5pt plus 0.3fil minus 2pt

{\bf Keywords:}\quad \thekeywords

\vglue 10pt plus 0.5fil minus 5pt

{\small  Proposed: \theproposer\hfill Received: \receiveddate\nl
Seconded: \theseconders\hfill 
\ifx\reviseddate\relax                         
Accepted: \accepteddate                        
\else
Revised: \reviseddate                          
\fi}
\eject
}       

\let\maketitlepage\maketitlep
\let\maketitle\maketitlepage

\font\phead=cmsl9 scaled 950
\font=cmsl9 scaled 1050
\font\pnum=cmbx10 scaled 913
\font\lnum=cmbx10 
\font\pfoot=cmsl9 scaled 950
\font=cmsl9 scaled 1050
\ifplaintex
\headline{\vbox to 0pt{\vskip -4.5mm\line{\small\phead\ifnum
\count0=\startpage ISSN 1364-0380 (on line)
1465-3060 (printed) \hfill {\pnum\folio}\else\ifodd\count0\def\\{ }%
\ifx\theshorttitle\relax\thetitle\else\theshorttitle\fi\hfill{\pnum\folio}
\else\def\\{ and }{\pnum\folio}\hfill\ifx\theshortauthors\relax\theauthors
\else\theshortauthors\fi\fi\fi}\vss}}
\footline{\vbox to 0pt{\vglue 0mm\line{\small\pfoot\ifnum\count0=\startpage
\copyright\ \gtp\hfill\else
\gt, Volume \thevolumenumber\ (\thevolumeyear)\hfill\fi}\vss
}}
\else
\makeatletter
\def\@oddhead{{\small\lhead\ifnum\count0=\startpage ISSN 1364-0380 (on line)
1465-3060 (printed) \hfill {\lnum\number\count0}\else\ifodd\count0
\def\\{ }\ifx\theshorttitle\relax \thetitle \else\theshorttitle\fi\hfill
{\lnum\number\count0}\else\def\\{ and }{\lnum\number\count0}
\hfill\ifx\theshortauthors\relax 
\theauthors\else\theshortauthors\fi\fi\fi}}\def\@evenhead{\@oddhead}
\def\@oddfoot{\small\lfoot\ifnum\count0=\startpage\copyright\ \gtp\hfill\else
\gt, Volume \thevolumenumber\ (\thevolumeyear)\hfill\fi}
\def\@evenfoot{\@oddfoot}
\makeatother
\fi

\newwrite\gtoutfile
\long\gdef\makeheadfile{  
{\def\\{, }\def\s{ }
\immediate\openout\gtoutfile head.xxx
\immediate\write\gtoutfile{Proxy-for: \ifx\theasciiauthors\relax
\theauthors\else\theasciiauthors\fi\s<\ifx\theasciiemail\relax\theemail\else\theasciiemail\fi>}
\immediate\write\gtoutfile{\noexpand\\}
\immediate\write\gtoutfile{Authors: \ifx\theasciiauthors\relax
\theauthors\else\theasciiauthors\fi}
{\def\\{ }\immediate\write\gtoutfile{Title: \ifx\theasciititle\relax
\thetitle\else\theasciititle\fi}}
\immediate\write\gtoutfile{Subj-class: GT or SG or MG etc}
\immediate\write\gtoutfile{MSC-class: \theprimaryclass\ifx\thesecondaryclass\relax\else, \thesecondaryclass\fi}
\immediate\write\gtoutfile{Journal-ref: Geom. Topol. \thevolumenumber
(\thevolumeyear) \startpage-\finishpage}
\immediate\write\gtoutfile{Comments: Published by Geometry and Topology at}
\immediate\write\gtoutfile{\s\s http://www.maths.warwick.ac.uk/gt/GTVol\thevolumenumber/paper\thepapernumber.abs.html}
\immediate\write\gtoutfile{\noexpand\\}
\immediate\write\gtoutfile{}
\ifx\theasciiabstract\relax
\immediate\write\gtoutfile{\theabstract}\else
\immediate\write\gtoutfile{\theasciiabstract}\fi
\immediate\write\gtoutfile{}
\immediate\write\gtoutfile{\noexpand\\}
\immediate\write\gtoutfile{}
\immediate\closeout\gtoutfile}}  

\def\maketitlepage{\maketitlep\makeheadfile}
\let\maketitle\maketitlepage

\lognumber{476}

\received{14 July 2004}
\volumenumber{9}\papernumber{22}\volumeyear{2005}
\pagenumbers{971}{990}   
\revised{18 May 2005}
\published{29 May 2005}
\accepted{18 May 2005}
\proposed{David Gabai}
\seconded{Walter Neumann, Cameron Gordon}

\usepackage{amsmath, amssymb}

\newcommand{\bd}{\partial}
\newcommand{\hm}{ho\-me\-o\-mor\-phic}
\newcommand{\R}{\ensuremath{\mathbb{R}}}
\newcommand{\RR}{\ensuremath{\mathbb{R}^2}}
\newcommand{\RRR}{\ensuremath{\mathbb{R}^3}}
\newcommand{\sbs}{\subseteq}
\newcommand{\irr}{ir\-re\-du\-ci\-ble}
\newcommand{\inc}{in\-com\-press\-i\-ble}

\newcommand{\ei}{end \irr}
\newcommand{\eei}{eventually \ei}
\newcommand{\tm}{3--manifold}

\newcommand{\p}{^{\prime}}
\newcommand{\pp}{^{\prime\prime}}
\newcommand{\er}{end reduction}
\newcommand{\ra}{\rightarrow}
\newcommand{\ns}{\emptyset}
\newcommand{\inte}{\operatorname{int}}
\newcommand{\ga}{\ensuremath{\gamma}}
\newcommand{\be}{\ensuremath{\beta}}
\newcommand{\al}{\ensuremath{\alpha}}

\newcommand{\n}{^{-1}}
\newcommand{\Inte}{\operatorname{Int}}

\newcommand{\de}{\ensuremath{\delta}}

\newcommand{\De}{\ensuremath{\Delta}}
\newcommand{\Dep}{\ensuremath{\De\p}}

\newcommand{\qe}{quasi-exhaustion}

\newcommand{\EE}{\ensuremath{\mathcal{E}}}
\newcommand{\BB}{\ensuremath{\mathcal{B}}}

\newcommand{\CC}{\ensuremath{\mathcal{C}}}

\newcommand{\whE}{\ensuremath{\widehat{E}}}
\newcommand{\whBB}{\ensuremath{\widehat{\BB}}}

\newcommand{\whB}{\ensuremath{\widehat{B}}}
\newcommand{\DD}{\ensuremath{\mathcal{D}}}
\newcommand{\whDD}{\ensuremath{\widehat{\DD}}}
\newcommand{\whSigma}{\ensuremath{\widehat{\Sigma}}}

\newcommand{\shp}{^{\#}}
\newcommand{\wh}{\widehat}
\newcommand{\wt}{\widetilde}
\newcommand{\wtV}{\ensuremath{\wt{V}}}
\newcommand{\wtM}{\ensuremath{\wt{M}}}
\newcommand{\wtDD}{\ensuremath{\wt{\DD}}}
\newcommand{\wtCC}{\ensuremath{\wt{\CC}}}
\newcommand{\wtA}{\ensuremath{\wt{A}}}
\newcommand{\wtD}{\ensuremath{\wt{D}}}
\newcommand{\wtB}{\ensuremath{\wt{B}}}
\newcommand{\wtBB}{\ensuremath{\wt{\BB}}}
\newcommand{\wtSigma}{\ensuremath{\wt{\Sigma}}}
\newcommand{\wttau}{\ensuremath{\wt{\tau}}}

\newcommand{\wtE}{\ensuremath{\wt{E}}}
\newcommand{\whk}{\ensuremath{\wh{k}}}
\newcommand{\whtau}{\ensuremath{\wh{\tau}}}

\newtheorem{thm}{Theorem}[section]

\newtheorem{lem}[thm]{Lemma}
\newtheorem{conj}[thm]{Conjecture}
\newtheorem{quest}[thm]{Question}

\begin{document}
\title[End reductions, fundamental groups and covering spaces]
{End reductions, fundamental groups, and covering\\spaces   
of irreducible open 3--manifolds}
\asciititle{End reductions, fundamental groups, and covering spaces   
of irreducible open 3-manifolds}
\author{Robert Myers}

\address{Department of Mathematics, Oklahoma State University\\Stillwater, 
OK 74078, USA}
\email{myersr@math.okstate.edu}
\urladdr{http://www.math.okstate.edu/~myersr/}

\begin{abstract} 
Suppose $M$ is a connected, open, orientable, irreducible 3--manifold which 
is not homeomorphic to \RRR. Given a compact 3--manifold $J$ in $M$ which 
satisfies certain conditions, Brin and Thickstun have associated to 
it an open neighborhood $V$  
called an \emph{end reduction} of $M$ at $J$. It has some useful
properties which allow one to extend to $M$ various results known to
hold for the more restrictive class of eventually end irreducible open
3--manifolds.

In this paper we explore the relationship of $V$ and $M$ with regard to their 
fundamental groups and their covering spaces. In particular we give 
conditions under which the inclusion induced homomorphism on fundamental 
groups is an isomorphism. We also show that if $M$  
has universal covering space homeomorphic to \RRR, then so does $V$. 

This work was motivated by a conjecture of Freedman (later disproved  
by Freedman and Gabai) on knots in $M$ which 
are covered by a standard set of lines in \RRR.  

\end{abstract}
\asciiabstract{%
Suppose M is a connected, open, orientable, irreducible 3-manifold
which is not homeomorphic to R^3.  Given a compact 3-manifold J in M
which satisfies certain conditions, Brin and Thickstun have associated
to it an open neighborhood V$ called an end reduction of M at J.  It
has some useful properties which allow one to extend to M various
results known to hold for the more restrictive class of eventually end
irreducible open 3-manifolds.  In this paper we explore the
relationship of V and M with regard to their fundamental groups and
their covering spaces.  In particular we give conditions under which
the inclusion induced homomorphism on fundamental groups is an
isomorphism.  We also show that if M has universal covering space
homeomorphic to R^3, then so does V.  This work was motivated by a
conjecture of Freedman (later disproved by Freedman and Gabai) on
knots in M which are covered by a standard set of lines in R^3.  }

\asciikeywords{3-manifold, end reduction, covering space}
\maketitlepage

\section{Introduction}

The Marden Conjecture states that an open hyperbolic 3--manifold with 
finitely generated fundamental group is almost compact, ie is 
homeomorphic to the interior of a compact 3--manifold. Independent 
proofs of this conjecture have recently been given by 
Ian Agol \cite{Ag} and by Danny Calegari and David Gabai \cite{CG}. 
These proofs use techniques which are mostly geometric.  
In \cite{Fr} Mike Freedman proposed the 
following topological conjecture which he proved implies the Marden 
Conjecture.  

\begin{conj}[Freedman Conjecture] Let $M$ be a connected, orientable 
open 3--manifold and $p\co\wtM\rightarrow M$ its universal covering map. 
Let \ga\ be a knot in $M$. If 
\begin{enumerate}
\item[{\rm(1)}] $\pi_1(M)$ is finitely generated, 
\item[{\rm(2)}] $\wtM$ is homeomorphic to \RRR, and 
\item[{\rm(3)}] $p\n(\ga)$ is a standard set of lines, 
\end{enumerate}
then $\pi_1(M-\ga)$ is finitely generated. \end{conj}

By a \textit{standard set of lines} we mean a subset $L$ of \RRR\ such 
that $(\RRR,L)$ is homeomorphic to $(\RR,X)\times\R$, where $X$ is 
a countably infinite closed discrete subset of \RR. 

In fact Freedman showed that the following special 
case of his conjecture implies the Marden Conjecture. 

\begin{conj}[Special Freedman Conjecture] Let $M$ and \ga\ be as above. 
Assume that, in addition,
\begin{enumerate}
\item[{\rm(4)}] $\pi_1(M)$ is a free product of two non-trivial groups, 
\item[{\rm(5)}] $[\ga]\in\pi_1(M)$ is algebraically disk busting, and 
\item[{\rm(6)}] $\langle\ga\rangle=0\in H_1(M,\mathbb{Z}_2)$. 
\end{enumerate}
Then $\pi_1(M-\ga)$ is finitely generated. \end{conj}

An element $g$ of a group $G$ is \textit{algebraically disk busting} if 
it is not conjugate into a proper free factor in any free factorization 
of $G$. For example, let $G$ be the free group on the set $\{a,b\}$, and 
let $g=a^2b^2$. If $g$ were conjugate into a proper free factor of $G$, 
then the quotient of $G$ by the normal closure of $g$  would have the form 
$\mathbb{Z}*\mathbb{Z}_n$, which contradicts the fact that 
this quotient is the fundamental group of the Klein bottle. 

This paper was originally part of a program for attacking 
Conjecture 1.2. After it was written Freedman and Gabai discovered 
a counterexample \cite{FG} to this conjecture, and so that program 
is now defunct. (We note that they had earlier found a counterexample 
to Conjecture 1.1.)  

The idea of the program was to split the problem into two parts: 
proving the conjecture 
for the case in which $M$ is ``end irreducible rel \ga'', and 
reducing the conjecture to this special case. (See the next 
section for definitions.)  This paper is concerned with the reduction 
to the special case. The Freedman--Gabai counterexample is 
end irreducible rel \ga, and so is a counterexample to the special case. 

The device for trying to do the reduction is an ``end reduction 
of $M$ at \ga''. This is a certain connected open subset $V$ of $M$ which 
is end irreducible rel \ga\ and has certain other nice properties. (The 
theory of end reductions was developed in more general contexts by Brin 
and Thickstun \cite{BT,BT2}, who were in turn inspired by earlier work 
of Brown and Feustel \cite{BF}.)  

The inclusion induced map $\pi_1(V)
\rightarrow\pi_1(M)$ is injective. Assuming the context of the 
Special Freedman Conjecture we prove in Theorem 5.1 that this 
homomorphism is also surjective, so $p\n(V)$ is connected;  
in Theorem 6.1 we prove that it is homeomorphic to \RRR\ and in Theorem 
8.1 that $p\n(\ga)$ is a standard set of lines in $p\n(V)$. 

The reduction would have worked as follows. 
Let $V$ be an end reduction of $M$ at \ga. Since $\pi_1(V)$ is isomorphic 
to $\pi_1(M)$ it is finitely generated and \ga\ is algebraically disk 
busting in $\pi_1(V)$. Then $p\n(V)$ is homeomorphic to \RRR\ and 
$p\n(\ga)$ is a standard set of lines in $p\n(V)$. If the conjecture 
were true for $V$ and \ga, then $\pi_1(V-\ga)$ would be finitely generated. 
The reduction would have been complete if the following has an 
affirmative answer. 

\begin{quest}Does $\pi_1(V-\ga)$ finitely generated imply that 
$\pi_1(M-\ga)$ is finitely generated? \end{quest} 

The paper is organized as follows. Sections 2 and 3 give an exposition of  
portions of Brin and Thickstun's theory of end reductions. 
Section 4 proves the incidental 
result that if $\pi_1(M)$ is finitely generated and indecomposable, then 
$V$ is either simply connected or $\pi_1$--surjective. Section 5 proves 
that if $\pi_1(M)$ is finitely generated and decomposable and $V$ is an 
end reduction of $M$ at an algebraically disk busting knot \ga, then 
$V$ is $\pi_1$--surjective. Section 6 proves that if $V$ is any end 
reduction of any irreducible, orientable  open 3--manifold $M$ whose 
universal covering 
space is homeomorphic to \RRR, then each component of $p\n(V)$ is 
homeomorphic to \RRR. Section 7 proves some lemmas about trivial 
$k$ component tangles which are then used in section 8 to prove that 
if $p\n(\ga)$ is a standard set of lines in $\wtM\approx\RRR$, then 
the intersection of each component of $p\n(V)$ with  
$p\n(\ga)$ is a standard set of lines.

Agol noticed that the proof that $\pi_1(V)$ and $\pi_1(M)$ are isomorphic 
when \ga\ is algebraically disk busting 
works in a more general context and used this in his proof of the 
Marden Conjecture. Section 9 shows how to modify the proof of Theorem 5.1 
to obtain a result which includes the situation considered by Agol. 
We note that Calegari and Gabai also used end reductions in their proof.      

This research was partially supported by NSF Grant DMS-0072429.

\section{End reductions}

In general we follow \cite{He} or \cite{Ja} for basic \tm\ terminology. 
When $X$ is a submanifold of $Y$ we denote the topological interior of 
$X$ by $\Inte X$ and the manifold interior of $X$ by $\inte X$. 
We say that $X$ is \textit{proper} in $Y$ if $X\cap C$ is compact 
for each compact $C\subseteq Y$ and that it is \textit{properly 
embedded} in $Y$ if $X\cap\bd Y=\bd X$. 
Two compact properly embedded surfaces $F$ and $G$ in a 3--manifold are in 
\textit{minimal general position} if they are in general position 
and for all surfaces isotopic to $F$ and in general position with 
respect to $G$ the intersection of $F$ and $G$ has the smallest 
number of components. 

Throughout the paper $M$ will be a connected, orientable, \irr, 
open \tm\ which is not \hm\ to \RRR.  A sequence $\{C_n\}_{n\geq0}$ 
of compact, connected \tm s $C_n$ in $M$ such that 
$C_n\sbs \inte C_{n+1}$ and $M-\inte C_n$ has no 
compact components is called a \textit{quasi-exhaustion in} $M$. 
If $\cup C_n=M$, then it is called an \textit{exhaustion for} $M$. 

In this section and the next we presents some basic material about 
Brin and Thickstun's theory of end reductions. See \cite{BT, BT2} for 
more details. We note that they work in a more general context and 
give somewhat different definitions of the terms that follow. In our 
context the definition of end reduction is equivalent to theirs.  

A compact, connected 3--manifold $J$ in $M$ is 
\textit{regular in} $M$ if $M-J$ is irreducible and has no component with 
compact closure. Since $M$ is \irr\ the first condition is equivalent to the 
statement that $J$ does not lie in a 3--ball in $M$. A quasi-exhaustion 
$\{C_n\}$ in $M$ is \textit{regular} if each $C_n$ is regular in $M$. 

Let $J$ be a regular 3--manifold in $M$, and let $V$ be an open subset 
of $M$ which contains $J$. We say that $V$ is \textit{\ei\ rel $J$ in 
$M$} if there is a regular \qe\ $\{C_n\}$ in $M$ such that 
$V=\cup_{n\geq0}C_n$, $J=C_0$, and $\bd C_n$ is \inc\ in $M-\inte J$ for all 
$n>0$. In the case that $V=M$ we say that $M$ is \textit{\ei\ rel} $J$;  
we say that $M$ is \textit{\eei\ } if it is \ei\ rel $J$ for some $J$.  

$V$ has the \textit{engulfing property rel $J$ in $M$} 
if whenever $N$ is regular in $M$, $J\sbs\inte N$, and $\bd N$ is \inc\ in 
$M-J$, then $V$ is ambient isotopic rel $J$ to $V\p$ such that $N\sbs V\p$. 
(Brin and Thickstun's ``weak engulfing'' property requires the isotopy 
to be fixed 
off a compact subset of $M-J$, but this can be achieved by using the 
covering isotopy theorem \cite{Ch, EK}.) 

$V$ is an \textit{\er\ of $M$ at $J$} if $V$ is \ei\ rel $J$ in $M$, $V$ 
has the engulfing property rel $J$ in $M$, and no component of $M-V$ is 
compact. (Note that \RRR\ contains no regular \tm s and hence has no 
end reductions, so we will henceforth not mention the hypothesis that 
$M$ is not \hm\ to \RRR.)

\begin{thm}[Brin--Thickstun] Given  a regular 3--manifold $J$ in $M$, 
an \er\ $V$ of $M$ at $J$ exists and is unique up to non-ambient 
isotopy rel $J$ in $M$. \end{thm}

\begin{proof} This follows from Theorems 2.1 and 2.3 of \cite{BT}. \end{proof}

We will need to understand the construction of $V$. 
Here is a brief sketch. We begin with a regular exhaustion 
$\{C_n\}_{n\geq0}$ of $M$ 
with $C_0=J$. Set $C_0^*=C_0$. If $\bd C_1$ is \inc\ in $M-J$ set $C_1^*=C_1$. 
Otherwise we ``completely compress'' $\bd C_1$ in $M-C_0^*$ to obtain $C_1^*$. 
This is done by  constructing a sequence of compact 3--manifolds 
$K_0, \ldots, K_p$ with $K_0=C_1$ and $K_p=C_1^*$ as follows. 
$K_{i+1}$ is obtained from $K_i$ by one of the following two operations.  
One may cut off a 1--handle from $K_i$ which misses $C_0$ and whose 
co-core meets $\bd K_i$ in an essential simple closed curve. 
(This is just compressing $\bd K_i$ ``to the inside''.) One may add 
a 2--handle to $K_i$ in $M-C_0$ whose core meets $\bd K_i$ in an essential 
simple closed curve. (This is just compressing $\bd K_i$ ``to the outside''.) 
Moreover the cutting and attaching occurs along 
disjoint annuli in $\bd C_1$. 
We may assume that $C_1^*\sbs\inte C_2$. If $\bd C_2$ is \inc\ in $M-J$ we 
set $C_2^*=C_2$. Otherwise we completely compress $\bd C_2$ in $M-C_1^*$ to 
get $C_2^*$. We continue in this fashion to construct a sequence 
$\{C_n^*\}_{n\geq0}$. We let $V^*=\cup_{n\geq0}C_n^*$ and then let $V$ be the 
component of $V^*$ containing $J$. $V^*$ is called a    
\textit{constructed end reduction of $M$ at $J$}. We will call $V$ a  
\textit{standard end reduction of $M$ at $J$}. 

It will be convenient to arrange for $V$ to be obtained solely through 
cutting 1--handles. The following result is in Brown and Feustel 
\cite{BF} and has antecedents in work of McMillan \cite{Mc, Mc2, Mc3}. 

\begin{lem}[Thick Exhaustion Lemma (Brown--Feustel)] 
There is an exhaustion $\{\wh{C}_n\}$ 
for $M$ such that $C_n\sbs \wh{C}_n$ and $C^*_n$ can be obtained from 
$\wh{C}_n$ by cutting 1--handles. \end{lem}

\begin{proof} Suppose that in the process of completely compressing 
$\bd C_n$ in $M-C_{n-1}$ we first cut some 1--handles and then add some 
2--handles. We may assume that the cores of these handles are in general 
position. If the intersection is non-empty, then the 2--handles cut the 
1--handles into a collection of shorter 1--handles. Thus one could obtain 
the same 3--manifold by first adding the 2--handles and then cutting a 
possibly greater number of 1--handles. 

It follows that we can add all the 2--handles at once to obtain a 
3--manifold $\wh{C}_n$ from which we then cut 1--handles to obtain $C^*_n$. 
\end{proof}

From now on we assume that our standard end reduction $V$ has been 
obtained by cutting 1--handles. 

It is clear that a standard end reduction $V$ is end irreducible 
rel $J$ in $M$. 

Suppose $N$ is regular in $M$, $J\sbs\inte N$, and $\bd N$ is 
incompressible in $M-J$. Then $N\sbs C_n$ for some $n$. We isotop 
$\bd N$ off the co-cores of the 1--handles which are cut to obtain 
$C_n^*$. Hence $N$ is now contained in the component of $C^*_n$ which 
contains $J$ and hence is contained in $V$. Running the isotopy backwards 
establishes the engulfing property. 

Suppose we give our cut 1--handles product structures of the form 
$D\times [a,b]$. We may arrange things so that if $H$ and $H\p$ are 
1--handles cut from $C_n$ and $C_m$, respectively, with $m>n$, then 
$H\cap H\p$ is either empty or consists of a finite number of 1--handles 
each of which has induced product structures from $H$ and from $H\p$. 
From this it follows that if $x\in C_n\cap(M-V^*)$, then $x$ lies in 
a properly embedded disk in $C_n$ which is in $M-V^*$. Also, if 
$x\in C_n\cap(V^*-V)$, then $x$ lies in a connected 3--manifold in 
$C_n$ which meets $\partial C_n$. It follows that $M-V$ has no compact 
component. 

Another consequence of the 1--handles having compatible product structures 
is the widely noted but seemingly unrecorded observation, 
first made by Brin and Thickstun, that $M-V^*$ is a lamination of $M$ 
each leaf of which is a plane. 

Now suppose that $J$ and $K$ are regular 3--manifolds in $M$ such that  
$J\sbs\inte K$ and $\bd K$ is incompressible in $M-J$. Let 
$V$ be an end reduction of $M$ at $J$ such that $K\sbs V$. 
In Theorem 2.2 of \cite{BT} Brin and Thickstun prove that if $L$ is 
regular in $M$, $K\sbs\inte L$, and $\bd L$ is incompressible in $M-K$, 
then $V$ is ambient isotopic rel $K$ to $V\p$ such that $L\sbs V\p$. 
This is called the \textit{strong engulfing property}. It follows from 
this property that $V$ is an end reduction of $M$ at $K$; see 
Corollary 2.2.1 of \cite{BT}. 

One can now prove   
the uniqueness of $V$ up to non-ambient isotopy rel $J$  
by using the strong engulfing property to inductively isotop 
the elements of properly chosen exhaustions $\{V_n\}$ for $V$ and 
$\{V_n\p\}$ for another end reduction $V\p$ of $M$ at $J$ so that 
$V_n\p\sbs V_n\sbs V_{n+1}\p$ for all $n$. See Theorem 2.3 of \cite{BT}. 

Suppose $V$ is an end reduction of $M$ at $J$ and $h_t\co V\rightarrow M$, 
$t\in[0,1]$,  
is a non-ambient isotopy with $h_0$ the inclusion map and $h_t(x)=x$ for 
all $x\in J$. By Lemma 2.5 of \cite{BT} $h_1(V)$ is an end reduction of 
$M$ at $J$. It is clear that $h_1(V)$ is end irreducible rel $J$ in $M$ 
and that $h_1(V)$ has the engulfing property rel $J$ in $M$. To see that 
$M-h_1(V)$ has no compact component we proceed as follows. 
(The author thanks Tom Thickstun for providing the following argument.) 
Let $\widehat{M}$ be the Freudenthal compactification of $M$. (See 
\cite{Ft} or section 1.3 of \cite{BT}.) The points of $\widehat{M}-M$ 
are the ends $E(M)$  of $M$. Since $M-V$ has no compact components each 
component 
of $\widehat{M}-V$ contains an end of $M$. It suffices to show that 
each component of $\widehat{M}-h_1(V)$ contains an end of $M$. 
Let $\{V_n\}$ be an exhaustion for $V$ such that each $V_n$ is regular 
in $M$. Thus each component of $\widehat{M}-V_n$ contains an end of $M$. 
By the covering isotopy theorem \cite{Ch,EK} there is an ambient isotopy 
$g_{n,t}\co\widehat{M}\rightarrow\widehat{M}$ rel $\widehat{M}-M$ such that 
$g_{n,t}(x)=h_t(x)$ for all $x\in V_n$. Therefore each component of 
$\widehat{M}-h_1(V_n)$ contains an end of $M$. Let $X$ be a component of 
$\widehat{M}-h_1(V)$. Then $X=\cap X_i$, where $X_i$ is a component of 
$\widehat{M}-h_1(V_i)$. Thus $X_i$ contains some $e_i\in E(M)$. We may 
assume that $e_i\rightarrow e\in \widehat{M}$. Since $E(M)$ is closed in 
$\widehat{M}$ and $X$ is the nested intersection of the $X_i$ we have 
that $e\in E(M)\cap X$. Thus $M-h_1(V)$ has no compact component.   

We now define a \textit{regular knot} in $M$ to be a simple closed 
curve \ga\ in $M$ which is not contained in a 3--ball in $M$. In this case 
a regular neighborhood $J$ of \ga\ in $M$ is a regular 3--manifold in $M$, 
and we define an \textit{end reduction of $M$ at \ga\ } to be an 
end reduction of $M$ at $J$. There is no loss of generality in working 
with end reductions at regular knots since every 
end reduction of $M$ at $J$ is also an end reduction at a regular knot 
in $M$. (See \cite{My Endcov}.)

\section[pi1-injectivity]{$\pi_1$--injectivity}

The following result is in Proposition 1.4 of \cite{BT2}. 
We include a proof here to make our exposition more self contained. 

\begin{thm}[Brin--Thickstun] Let $M$ be a connected, orientable, \irr, 
open \tm. Let $J$ 
be a regular 3--manifold in $M$ and $V$ an end reduction of $M$ at $J$. 
Then the inclusion induced homomorphism $\pi_1(V)\rightarrow\pi_1(M)$ 
is injective. \end{thm}

\begin{proof} Since $V$ is unique up to non-ambient isotopy rel $J$ 
we may assume that $V$ is a standard end reduction associated to 
an exhaustion $\{C_n\}$ for $M$. By Lemma 2.2 we may assume that 
$C_n^*$ is obtained from $C_n$ by cutting 1--handles. 

Now suppose that \al\ is a loop in $V$ which bounds a singular disk 
in $M$. Thus we have a map $f\co\De\rightarrow M$, where \De\ is a disk 
and $f(\bd\De)=\al$. 

We may assume that $\al\sbs C^*_1$ and $f(\De)\sbs C_1$. 
Let $D$ be the co-core of the first 1--handle which is cut from $C_1$ in the 
process of obtaining $C^*_1$. Put $f$  in general position with respect 
to $D$. Then $f\n(D)$ consists of simple closed curves in $\inte\De$. 
We redefine $f$ in a neighborhood of the outermost disks among these curves 
so that $f\n(D)=\ns$. 

We then apply this procedure to the next 1--handle which is cut and 
continue until we have that $f(\De)\sbs C^*_1$ and hence $f(\De)\sbs V$. 
\end{proof}

\section[pi1-surjectivity: the indecomposable case]{$\pi_1$--surjectivity: the indecomposable case}

In general, end reductions need not be $\pi_1$--surjective. In this section 
and the next  
we restrict to the case that $\pi_1(M)$ is finitely generated and place 
further restrictions which guarantee $\pi_1$--surjectivity. 

Recall that the Scott compact core theorem \cite{Sc} asserts that if 
$M$ is a connected non-compact 3--manifold with $\pi_1(M)$ finitely generated, 
then there is a compact, connected 3--manifold $N$ in $\inte M$ 
(a \textit{compact core}) such that the inclusion induced homomorphism 
$\pi_1(N)\rightarrow\pi_1(M)$ is an isomorphism. Since we will be  assuming 
that $M$ is \irr\ and $\pi_1(M)$ is non-trivial we may cap off any 
2--sphere components of $\bd N$ with 3--balls in $M-\inte N$ and hence may 
assume that $\bd N$ contains no 2--spheres. 

Note that if we are in a situation in which we wish to show that 
$\pi_1(V)\rightarrow\pi_1(M)$ is an isomorphism it suffices to show 
that there is a compact core $N$ such that $N\sbs V$. 

\begin{thm}  Let $M$ be a connected, orientable, \irr, open \tm. Let $J$ 
be a regular 3--manifold in $M$ and $V$ an end reduction of $M$ at $J$. 
Suppose that $\pi_1(V)$ is non-trivial and that $\pi_1(M)$ is finitely 
generated, non-cyclic, and indecomposable. Then the inclusion induced 
homomorphism 
$\pi_1(V)\rightarrow\pi_1(M)$ is an isomorphism. \end{thm}
 
\begin{proof} We may assume that $V$ is a standard end reduction 
obtained from an exhaustion $\{C_n\}$ for $M$. We may further assume 
that $C_n^*$ is obtained from $C_n$ by cutting 1--handles and that 
a compact core $N$ lies in $\inte C_1$. 
Since $\pi_1(M)$ is indecomposable and non-cyclic we have that $\bd N$ 
is \inc\ in $M$. 

Let $D$ be the co-core of the first 1--handle of $C_1$ which is cut. 
Isotop $N$ in $\inte C_1$ so that $\bd N$ is in general position with 
respect to $D$ 
and meets it in a minimal number of components. Each of these components 
bounds a disk on $\bd N$. Let \De\ be an innermost such disk. Then 
$\bd\De=\bd\Dep$ for a disk \Dep\ on $D$. The 2--sphere $\De\cup\Dep$ bounds 
a 3--ball $B$ in $M$. Since $C_1$ is regular it is \irr, and so $B$ lies in 
$\inte C_1$. One may then isotop \De\ across $B$ past \Dep\ to remove at 
least one component of $D\cap \bd N$, contradicting minimality. 
Therefore $D\cap \bd N=\ns$. 

In a similar fashion we see that $N$ meets none of the co-cores of the 
1--handles which are cut to get $C^*_1$. Thus $N$ lies in $C^*_1$ and 
hence in $V^*$. 

Suppose $N$ does not lie in $V$. Let \al\ be a loop in $V$. Then \al\ 
is freely homotopic in $M$ to a loop in $N$. We may assume that \al\ and 
the image of the homotopy lie in $C_2$. We may further assume that \al\ 
and $N$ lie in different components of $C_2^*$. It follows that the 
homotopy can be cut off on the co-core of some 1--handle which is cut 
to obtain $C^*_2$, and thus \al\ is homotopically trivial. Hence 
$\pi_1(V)$ is trivial, contradicting our hypothesis. 

Thus $N$ lies in $V$, which completes the proof. \end{proof}

\section[pi1-surjectivity: the algebraically disk busting case]{$\pi_1$--surjectivity: the algebraically disk busting case}

\begin{thm}  Let $M$ be a connected, orientable, \irr, open \tm. Let $\gamma$ 
be a regular knot in $M$ and $V$ an end reduction of $M$ at $\gamma$. 
Suppose that $\pi_1(M)$ is finitely generated and decomposable and that 
$\gamma$ is algebraically disk busting.  
Then the inclusion induced homomorphism 
$\pi_1(V)\rightarrow\pi_1(M)$ is an isomorphism. \end{thm}

\begin{proof} As before we may assume that $V$ is a standard end reduction 
obtained from the  
exhaustion $\{C_n\}$ by cutting 1--handles. Let $N$ be a compact core 
for $M$. We may assume that $N\sbs \inte C_1$ and that \ga\ is freely 
homotopic in $\inte C_1$ to a knot $\gamma\p\sbs\inte N$. 

Let $D$ be the co-core of the first 1--handle which is cut from $C_1$. 
Let $C_1\p$ be the resulting \tm. Note that $\gamma\cap D=\ns$. 

$N$ is either a cube with handles or consists of a finite number of 
3--manifolds with boundaries which are \inc\ in $M$  to which 1--handles 
have been added. It follows that $N$ can be chosen so that $N\cap D$ 
consists of essential disks in $N$. Assume that among all choices of $N$ 
the number of these disks is minimal. 

Suppose $N\cap D\neq\ns$. Isotop $\gamma\p$ in $\inte C_1$ so that 
$\gamma\p\sbs \inte N$, $\gamma\p$ is in general position with respect 
to $D$ and $\gamma\p\cap D$ has a minimal number of points. 

If $\gamma\p\cap D=\ns$, then some essential disk in $N$ misses $\gamma\p$. 
It follows that $\gamma\p$ is conjugate into a proper free factor of 
$\pi_1(N)$ and hence that $\gamma$ is not algebraically disk busting in 
$\pi_1(M)$, a contradiction. 

Thus $\gamma\p\cap D\neq\ns$. We will show that, contrary to our assumptions, 
there is a compact core $N\p$ in $C_1$ which meets $D$ in fewer disks 
than $N$ does. 

We have a map $f\co S^1\times[0,1]\rightarrow C_1$ with $f(S^1\times\{0\})=
\gamma$ and $f(S^1\times\{1\})=\gamma\p$. Put $f$ in general position with 
respect to $D$. Then $f\n(D)$ consists of simple closed curves and arcs with 
boundary in $S^1\times\{1\}$. Since  $\gamma\p\cap D\neq\ns$ there must be 
a component of the latter type. Its union with some arc on $S^1\times\{1\}$ 
bounds a disk in $S^1\times[0,1]$. We may assume that the interior of this 
disk is disjoint from $f\n(D)$. The disk provides a path homotopy in $C_1$ 
between an arc \al\ in $\ga\p$ and a path \be\ in $D$. 

\textbf{Case 1}\qua $\bd\be$ lies in a single component \De\ of $N\cap D$. 

Then \be\ is homotopic rel $\bd\be$ in $D$ to an arc \de\ in \De\ 
such that $\de\cap\ga\p=\bd\de$. Let $\ga\pp=\ga\p\cup\de-\al$. 
Isotop $\ga\pp$ slightly so that \de\ is moved off $D$ and 
$\ga\pp\cap D=(\gamma\p\cap D)-\bd\al$. Thus $\gamma\pp$ meets $D$ in 
fewer points than did $\gamma\p$, a contradiction. 

\textbf{Case 2}\qua $\bd\be$ lies in two different components $\De_0$ and 
$\De_1$ of $N\cap D$. 

Then \be\ is homotopic rel $\bd\be$ in $D$ to an arc \de\ in $D$ which 
meets each $\De_i$ in an arc $\de_i$ and is otherwise disjoint from $N$. 
Push \de\ slightly off $D$ and add a 1--handle to $N$ whose core is 
$\de\cap(C_1-\Inte N)$. Then cut a 1--handle from $N$ whose co-core is 
$\De_1$. Call the result $N\p$. Since \al\ is path homotopic to \be\ 
the core of the cut 1--handle is path homotopic to an arc in $N\p$, from 
which it follows that $N\p$ is a compact core for $M$. Although it 
no longer contains $\ga\p$ it meets $D$ in one fewer disk, contradicting 
our assumption on $N$. 

Thus we must have that $N\cap D=\ns$. We repeat the argument with 
$C\p_1$ in place of $C_1$ and continue until we get that $N\sbs C^*_1$. 
As in Theorem 4.1 we get that $N\sbs V$, and we are done. \end{proof}

\section[V is homeomorphic to R3]{$\wtV$ is \hm\ to \RRR}

\begin{thm} Let $M$ be a connected, orientable, \irr, open \tm. Let $J$ 
be a regular 3--manifold in $M$ and $V$ an end reduction of $M$ at $J$. 
Let $p\co\wtM\ra M$ be the universal covering map. Let \wtV\ be a 
component of $p\n(V)$. If \wtM\ is \hm\ to \RRR, then so is \wtV. 
\end{thm}
 
\begin{proof} Since $V$ is $\pi_1$--injective in $M$ we have that \wtV\ is 
simply connected. It suffices to show that for 
every compact, connected 
subset $A$ of \wtV\ there exists a 3--ball in \wtV\ whose interior 
contains $A$. 

We may assume that $V$ is the end reduction obtained from an exhaustion 
$\{C_n\}$ for $M$ by cutting 1--handles. 

Since \wtM\ is \hm\ to \RRR, there is a 3--ball $B$ in \wtM\ whose interior 
contains $A$. We may assume that $(p(B),p(A))\sbs(C_1,C_1^*\cap V)$. 

Suppose $D$ is the co-core of the first 1--handle of $C_1$ which is cut. 
So $D$ compresses $\bd C_1$ in $C_1-J$. Let $C_1\p$ be the resulting \tm.  
Then $p\n(D)$ is a disjoint union 
of disks $\wt{D}_i$ in $p\n(C_1)$. Since the group of covering translations 
is properly discontinuous only finitely many of these disks meet $B$. 
None of them meets $A$. Put $\bd B$ in general position with respect to 
$p\n(D)$.  

Assume that for some $i$ we have $\wt{D}_i\cap B\neq\ns$. Let \De\ be an   
innermost disk on $\wtD_i$, ie $\De\cap\bd B=\bd\De$. 

Suppose \De\ lies in $B$. Then it splits $B$ into two 3--balls one of which, 
say $B\shp$, contains $A$. Then $B\shp$ can be isotoped to remove $\bd\De$ 
from the intersection. 

Now suppose \De\ lies in $\wtM-\inte B$. Attaching a 2--handle to $B$ 
with core \De\ gives a \tm\ $Q$ with $\bd Q$ a disjoint union of two 
2--spheres. One of them bounds a 3--ball $B\shp$ in \wtM\ which contains $Q$. 
Then $\bd B\shp$ has at least one fewer intersection with $p\n(D)$ than 
did $\bd B$. 

Continuing in this fashion we get a 3--ball $B\p$ in \wtM\  such that 
$B\p\cap p\n(D)=\ns$. So $A\sbs\inte B\p$ and $B\p\sbs\inte p\n(C_1\p)$. 

We then proceed to the cutting of the next 1--handle. Eventually we get 
a 3--ball $B^*$ in \wtM\ such that $A\sbs \inte B^*$ and $B^*\sbs 
p\n(C_1^*\cap V)\sbs\wtV$. \end{proof}

\section{Standard embeddings of trivial $k$--tangles}

A \textit{$k$--tangle} is a pair $(B,\tau)$, where $B$ is a 3--ball and 
$\tau$ is a union of $k$ disjoint properly embedded 
arcs $\tau_i$ in $B$. $(B,\tau)$ and $(B\p,\tau\p)$ have the same 
\textit{tangle type} if there is a homeomorphism of pairs 
$h\co(B,\tau)\ra(B\p,\tau\p)$.  
We say that $(B,\tau)$ is \textit{trivial} 
if it has the same tangle type as $(D,P)\times[0,1]$, 
where $D$ is a disk and $P$ is a set of $k$ points in $\inte D$. 

A \textit{shell} is a \tm\ $\Sigma$ which is \hm\ to $S^2\times[0,1]$. 
Let $\beta$ be a union of $k$ disjoint properly embedded arcs $\beta_i$ in 
$\Sigma$ such that each $\beta_i$ joins $S^2\times\{0\}$ to $S^2\times\{1\}$.  
$(\Sigma,\beta)$ is a $k$--\textit{braid} if there is a homeomorphism 
$g\co(\Sigma,\beta)\ra(S^2,P)\times[0,1]$, where $P$ is a 
set of $k$ points in $S^2$. 

\begin{lem} 

{\rm(1)}\qua $(B,\tau)$ is a trivial $k$--tangle if and only if there 
is a disjoint union $F$ of $k$ disks $F_i$ in $B$ such that for 
$1\leq i\leq k$ one has that $F_i\cap\tau=\tau_i$ and 
$F_i\cap\bd B=\bd F_i-\inte\tau_i$. 

{\rm(2)}\qua $(\Sigma,\beta)$ is a $k$--braid if and only if for some 
(and hence any) $\beta_j$ its  
exterior is a 3--ball $X$ such that either $k=1$ or $k>1$ and $\beta-\beta_j$ 
is a trivial $(k-1)$--tangle in $X$ for which there exist disks $F_i$ as 
in (1) each of which has connected intersection with the annulus 
$\bd X-\inte(X\cap\bd\Sigma)$. \end{lem}

\begin{proof} The proof is left as an exercise. \end{proof}

\begin{lem} 

{\rm(1)}\qua Suppose $(B,\tau)$ is a trivial $k$--tangle and $G$ is a
finite disjoint union of properly embedded disks in $B$ which misses
$\tau$.  Let $B\p$ be the closure of a component of $B-G$ such that
$B\p\cap\tau\neq\ns$. Then $(B\p,B\p\cap\tau)$ is a trivial $k\p$--tangle
for some $k\p\leq k$.

{\rm(2)}\qua Suppose $(\Sigma,\beta)$ is a $k$--braid and $G$ is a
finite disjoint union of properly embedded disks in $S^2\times (0,1]$
which misses $\beta$.  Let $\Sigma\p$ be the closure of the component
of $\Sigma-G$ which contains $S^2\times\{0\}$. Then $(\Sigma\p,\beta)$
is a $k$--braid. \end{lem}

\begin{proof} (1)\qua Let $F$ be as in Lemma 7.1(1). Put $F$ in minimal 
general position with respect to $G$. Since the exterior of $\tau$ in $B$ 
is a cube with handles it is irreducible; it follows that $G\cap F$ has 
no simple closed curve components. Let $\tau\p=\tau\cap B\p$. 
For each component $\tau_i$ of $\tau\p$ let $F_i\p$ be the closure of 
the component of $F_i-(F_i\cap G)$ which contains $\tau_i$. It then 
follows from Lemma 7.1 (1) that $(B\p,\tau\p)$ is trivial. 

(2)\qua Let $\beta_j$, $X$, and $F$ be as in Lemma 7.1(2). 
Again putting $F$ in minimal general 
position with respect to $G$, no component of $G\cap F$ is a simple 
closed curve. We have that $\beta$ lies in $\Sigma\p$. For each component 
$\beta_i$ of $\beta-\beta_j$ let $F_i\p$ be the closure of the component 
of $F_i-(F_i\cap G)$ which contains $\beta_i$. Then $F_i\p$ lies in 
$X\p=X\cap\Sigma\p$ and has connected intersection with 
$\bd X\p-\inte(X\cap\bd\Sigma\p)$. Thus by Lemma 7.1 (2) $(\Sigma\p,\beta)$ 
is a $k$--braid. \end{proof}

Let $B$ and $\whB$ be 3--balls, $\tau$ a $k$--tangle in $B$, and 
$\whtau$ a $\whk$--tangle in $\whB$. We say that $(B,\tau)$ is 
\textit{standardly embedded} in $(\whB,\whtau)$ if 
(i) $B\sbs\inte\whB$, (ii) $\whtau\cap B=\tau$, and (iii) there is a finite 
disjoint union \whDD\ of properly embedded disks in $\whB-B$ such that 
$\whDD\cap\whtau=\ns$ and \whDD\ splits $\whB-\inte B$ into a shell 
$\whSigma$ which 
meets \whtau\ in a $2k$--braid and a disjoint union $\whBB$ of 3--balls 
$\whB_j$  
each of which meets \whtau\ in a trivial $k_j$--tangle ($k_j$ depending 
on the 3--ball $\whB_j$). 

\begin{lem} Suppose $(B,\tau)$ is standardly embedded in $(\whB,\whtau)$. 
Then 

{\rm(1)}\qua one may choose \whDD\ (and hence $\whBB$) to be connected,  

{\rm(2)}\qua for each component $\whtau_i$ of $\whtau$ one 
has that $\whtau_i\cap B$ is connected,  

{\rm(3)}\qua if $(B,\tau)$ is trivial then so is $(\whB,\whtau)$, and 

{\rm(4)}\qua if $(\whB,\whtau)$ is standardly embedded in $(\wtB,\wttau)$, 
then $(B,\tau)$ is standardly embedded in $(\wtB,\wttau)$. \end{lem}

\begin{proof} (1) and (3) follow from the fact that if 
$(B\p,\tau\p)$ is a trivial $k\p$--tangle and $(B\pp,\tau\pp)$ is a 
trivial $k\pp$--tangle such that $B\p\cap B\pp$ is a disk which is 
disjoint from $\tau\p\cup\tau\pp$, then 
$(B\p\cup B\pp,\tau\p\cup\tau\pp)$ is a 
trivial $(k\p+k\pp)$--tangle. 

(3) follows from the fact that 
$\whtau\cap\Sigma$ is a $2k$--braid. 

By (1) we may assume that \whDD\ is connected and that 
there is a properly embedded disk $\wtDD$ in $\wtB-\whB$ which 
splits $\wtB-\inte\whB$ into a shell $\wtSigma$ and a 3--ball $\wtBB$ 
such that $(\wtSigma,\wtSigma\cap\wttau)$ is a $2\wt{k}$--braid and 
$(\wtBB, \wtBB\cap\wttau)$ is a trivial $(\wt{k}-\wh{k})$--tangle. 
Let $\wtA=(\bd\whDD)\times[0,1]$ in the product structure on \wtSigma. 
Since \wttau\ misses \wtDD\ we may isotop \wtA\ rel $\bd\whDD$ so 
that it misses \wtDD. Then $\whDD\cup\wtA$ splits $\wtB-\inte B$ into 
a shell which meets \wttau\ in a $2k$--braid and a 3--ball which meets 
\wttau\ in a trivial $(\wt{k}-k)$--tangle. Thus we have (4). \end{proof} 

\begin{lem} Suppose $(B,\tau)$ is standardly embedded in $(\whB,\whtau)$. 
Let $\whE$ be a properly embedded disk in $\whB-B$ which is disjoint from 
\whtau. 
Let $B\shp$ be the closure of the component of $\whB-\whE$ which contains $B$. 
Let $\tau\shp=\whtau\cap B\shp$. Then $(B,\tau)$ is standardly embedded 
in $(B\shp,\tau\shp)$. \end{lem}

\begin{proof} Since the definition of standardly embedded does not depend 
on the tangle type of $(B,\tau)$ we may assume that $(B,\tau)$ is trivial. 
By Lemma 7.3 (1) we may assume that \whDD\ is connected. Put it in 
minimal general position with respect to $\whE$. The exterior of 
$\whtau\cap(\whB-\inte B)$ in $\whB-\inte B$ is a cube with handles 
and so is \irr. It follows that $\whDD\cap \whE$ has no simple closed curve 
components. 

Thus $\whDD\cap B\shp$ is a disjoint union $\DD\shp$ of disks. 
These disks split $B\shp-\inte B$ into a shell $\Sigma\shp$ and a disjoint union 
$\BB\shp$ of 3--balls. 
By deleting some of these disks we may assume that each of these 
3--balls meets $\tau\shp$. By Lemma 7.3 (3) we have that  $(\whB,\whtau)$ is 
trivial, so by Lemma 7.2 (1) $(B\shp,\tau\shp)$ is trivial. A second application 
of Lemma 7.2 (1) shows that the intersection of $\tau\shp$ with 
each component of $\BB\shp$ is a trivial tangle in that 3--ball. 
By Lemma 7.3 (2) we have that $(\Sigma\shp,\Sigma\shp\cap\tau\shp)$ 
is a $2k$--braid. Hence $(B,\tau)$ is standardly embedded in $(B\shp,\tau\shp)$. 
\end{proof}

\begin{lem} Let $(B,\tau)$, $(\whB,\whtau)$, and $(B\shp,\tau\shp)$ be 
as in the previous lemma. Suppose $(\whB,\whtau)$ is standardly embedded 
in $(\wt{B},\wt{\tau})$. Then $(B\shp,\tau\shp)$ is standardly embedded in 
$(\wt{B},\wt{\tau})$. \end{lem}

\begin{proof} Recall that $B\shp$ is the closure of the component of 
$\whB-\whE$ which contains $B$. Since $(\whB,\whtau)$ is standardly 
embedded in $(\wtB,\wttau)$ there is a properly embedded disk $\wtDD$ 
in $\wtB-\whB$ which splits 
$\wtB-\inte\whB$ into a shell $\wt{\Sigma}$ and a 3--ball \wtCC. 

Let $\wtA$ be the annulus $\bd\whE\times[0,1]$ in 
$\wt{\Sigma}=(\bd\whB)\times[0,1]$. Note that we may have 
$\wtA\cap\wtDD\neq\ns$. 
Since $\wttau\cap\wtDD=\ns$ we may isotop $\wtA$ in $\wt{\Sigma}-\wttau$ so 
that $\wtA\cap\wtDD=\ns$. We then deform the product structure on 
$\wt{\Sigma}$ so that we still have $\wtA=\bd\whE\times[0,1]$. 

Now let \wtE\ be a parallel copy of $\whE\cup \wtA$ in the component of 
$\wtB-(\whE\cup \wtA)$ which does not contain $B$. We may assume that 
$\wtE\cap(\wttau\cup\wtDD)=\ns$. Now \wtE\ splits $\wtB-\inte B\shp$ into a 
shell $\Sigma^+$ and a 3--ball $\CC^+$. 

Note that $\Sigma^+$ is the union of the parallelism 
$(\whE\cup \wtA)\times[0,1]$ 
between $\whE\cup \wtA$ and \wtE\ and $(\bd B\shp-\inte\whE)\times[0,1]$, where 
the latter product structure is that of $\wt{\Sigma}$. It follows that 
$(\Sigma^+,\Sigma^+\cap\wttau)$ is a braid. 

Let $\CC^+_0$ be the closure of the component of $\wtB-(\wtE\cup\wtDD)$ 
such that $\bd \CC^+_0$ contains $\wtE\cup\wtDD$. Then $\CC^+_0\cap\whB$ is 
a 3--ball which contains $\whtau-\tau\shp$. As in the proof of 
Lemma 7.4 we may assume that $(B,\tau)$ is trivial, and thus by Lemma 7.3 (3) 
so is $(\whB,\whtau)$, and thus by Lemma 7.2 (1) so is 
$(\CC^+_0\cap\whB,\whtau-\tau\shp)$. 
Now $\CC^+_0-\Inte\whB$ is a 3--ball whose intersection with \wttau\ consists 
of product arcs in the product structure induced by that of $\wt{\Sigma}$. 
It follows that $(\CC^+_0,\CC^+_0\cap\wttau)$ is a trivial tangle. Since 
$(\wtCC,\wtCC\cap\wttau)$ is trivial by the definition of standardly embedded, 
we have that $(\CC^+_0\cup\wtCC,\wttau-\tau\shp)$ is trivial. Hence 
$(B\shp,\tau\shp)$ is standardly embedded in $(\wtB,\wttau)$. \end{proof}

\section[V cap p-inv(g) is standard in V]{$\wt{V}\cap p\n(\ga)$ is standard in $\wt{V}$}

Recall that a proper 1--manifold $L$ in \RRR\ is a 
\textit{standard set of lines} if there is a homeomorphism 
$H\co(\RRR,L)\ra(\RR,X)\times\R$, where $X$ is a countably 
infinite closed discrete subset of \RR. 

\begin{thm} Let $M$ be a connected, orientable, irreducible open 
3--manifold. Let \ga\ be a regular knot in $M$ and $V$ an end reduction 
of $M$ at \ga. Let $p\co\wt{M}\ra M$ be the universal covering map. 
Suppose $\wt{M}$ is \hm\ to \RRR\ and $p\n(\ga)$ is a standard set 
of lines in $\wt{M}$. Then for each component $\wt{V}$ of $p\n(V)$ 
we have that $\wt{V}\cap p\n(\ga)$ is a standard set of lines in 
$\wt{V}$. \end{thm}

Recall that by Theorem 6.1 $\wt{V}$ is \hm\ to \RRR. 
We will make use of the following lemma. 

\begin{lem} $L$ is a standard set of lines in \RRR\ if and only if there 
is an  
exhaustion $\{B_n\}_{n\geq0}$ for \RRR\ with each $B_n$ a 3--ball such that  

{\rm(1)}\qua for each $n\geq0$ $B_n\cap L$ is a $k_n$--tangle $\tau^n$,  

{\rm(2)}\qua $(B_0,\tau^0)$ is trivial, and 

{\rm(3)}\qua for each $n\geq0$ $(B_n,\tau^n)$ is standardly embedded in 
$(B_{n+1},\tau^{n+1})$. \end{lem}

\proof
Let $q\co\RR\times\R\ra\RR$ be projection onto the first factor. 

First suppose that $L$ is standard. 

Choose an exhaustion $\{E_n\}_{n\geq0}$ of \RR\ such that 
each $E_n$ is a disk, 
$E_0\cap q(H(L)) \neq \ns$, $(\bd E_n)\cap q(H(L))=\ns$, and 
$(E_{n+1}-E_n)\cap q(H(L))\neq\ns$ for all $n\geq0$. Choose a properly 
embedded arc $\al_{n+1}$ in $E_{n+1}-E_n$ which splits $E_{n+1}-\inte E_n$ 
into a disk $\EE_{n+1}$ and an annulus $A_{n+1}$ such that 
$A_{n+1}\cap q(H(L))=\ns$. 

We define $B_n$, $\DD_{n+1}$, $\BB_{n+1}$, and $\Sigma_{n+1}$ as follows. 
\begin{align*}
B_n&=H^{-1}(E_n\times[-(n+1),n+1]) \\
\mathcal{D}_{n+1}&=H^{-1}(\al_{n+1}\times[-(n+2),n+2]) \\
\mathcal{B}_{n+1}&=H^{-1}(\EE_{n+1}\times[-(n+2),n+2]) \\
\Sigma_{n+1}&=B_{n+1}-(\inte B_n\cup \Inte \mathcal{B}_{n+1})
\end{align*}
Then the exhaustion $\{B_n\}_{n\geq0}$ of \RRR\ has 
the required properties. In particular the existence of $\DD_{n+1}$, 
$\Sigma_{n+1}$, and $\BB_{n+1}$ establishes (3). 

Now suppose that $\{B_n\}_{n\geq0}$ is an exhaustion satisfying the 
three properties. 

We construct $H$ inductively. 
Let $\{E_n\}_{n\geq0}$ be an exhaustion of \RR\ by disks. 

Since $B_0\cap L$ is a trivial $k_0$--tangle there is a 
homeomorphism $H_0\co(B_0,B_0\cap L)\ra(E_0,P_0)\times[-1,1]$, 
where $P_0$ is a set of $k_0$ points in $\inte E_0$. 

Suppose we have constructed a homeomorphism 
\[H_n\co(B_n,B_n\cap L)\ra(E_n,P_n)\times[-(n+1),n+1],\] 
where $P_n$ is a set of $k_n$ points in $\inte E_n$. 

Choose a properly embedded arc $\al_{n+1}$ in $E_{n+1}-E_n$. Then $\al_{n+1}$ 
splits $E_{n+1}-\inte E_n$ into an annulus $A_{n+1}$ and a disk 
$\EE_{n+1}$. Suppose $B_{n+1}\cap L$ has $k_{n+1}$ components. 
Let $P_{n+1}$ be the union of $P_n$ and a set of $k_{n+1}-k_n$ points 
in $\inte \EE_{n+1}$. 

By (3) and Lemma 7.3 (1) there is a properly embedded disk $\DD_{n+1}$ in 
$B_{n+1}-B_n$ which splits $B_{n+1}-\inte B_n$ into a shell $\Sigma_{n+1}$ 
and a 3--ball $\BB_{n+1}$ such that $(\Sigma_{n+1},\Sigma_{n+1}\cap L)$ is 
a $2k_n$--braid and $(\BB_{n+1},\BB_{n+1}\cap L)$ is a trivial 
$(k_{n+1}-k_n)$--tangle. 

We define a shell $S_{n+1}$ as follows. \[S_{n+1}=
((E_n\cup A_{n+1})\times[-(n+2),n+2)])-\inte(E_n\times[-(n+1),n+1])\] 
The restriction of $H_n$ to $(\bd B_n,\bd B_n\cap L)$ extends to a 
homeomorphism 
\[H^{\Sigma}_{n+1}\co(\Sigma_{n+1},\Sigma_{n+1}\cap L)\ra 
(S_{n+1},S_{n+1}\cap(P_n\times[-(n+2),n+2])).\]
The restriction of $H_{n+1}^{\Sigma}$ to \[\Sigma_{n+1}\cap\mathcal{B}_{n+1}=
(H^{\Sigma}_{n+1})^{-1}(\al_{n+1}\times[-(n+2),n+2])\] then extends to a 
homeomorphism 
\[H_{n+1}^{\BB}\co(\mathcal{B}_{n+1},\mathcal{B}_{n+1}\cap L)\ra
(\EE_{n+1},P_{n+1}-P_n)\times[-(n+2),n+2].\]
We then use $H_n$, $H_{n+1}^{\Sigma}$, and $H^{\BB}_{n+1}$ to define 
\[H_{n+1}\co(B_{n+1},B_{n+1}\cap
L)\ra(E_{n+1},P_{n+1})\times[-(n+2),n+2].\eqno{\qed}\] 

\begin{proof}[Proof of Theorem 8.1] Given an exhaustion 
$\{B_n\}$ for \wtM\ which satisfies the conditions of Lemma 8.2, we will 
construct an exhaustion $\{B_n^*\}$ for \wtV\ which satisfies these 
conditions. We may assume that $B_0\sbs\wtV$. Let $B_0^*=B_0$. 

The proof follows that of Theorem 6.1. We assume that $V$ is the end 
reduction obtained from an exhaustion $\{C_n\}$ for $M$ by cutting 
1--handles. 

Suppose $A_0$ is a compact, connected subset of \wtV. We may assume that 
the interior of $B_1$ contains $A_0\cup B_0^*$ and that 
$(p(B_1),p(A_0\cup B_0^*))\sbs (C_1,C_1^*\cap V)$. 

Let $D$ be the co-core of the first 1--handle of $C_1$ which is cut.  
Then $D$ compresses $\bd C_1$ in $C_1-\gamma$. Let $C_1\p$ be the 
resulting 3--manifold. Then $p\n(D)$ is a disjoint union of disks 
$\wt{D}_i$ in $p\n(C_1)$. Only finitely many $\wtD_i$ meet $B_1$. 
None of them meet $A_0\cup B_0^*\cup p\n(\ga)$. Put $\bd B_1$ in general 
position with respect to $p\n(D)$. Assume that for some $i$ we have 
$\wtD_i\cap B\neq\ns$. Let \De\ be an innermost disk on $\wtD_i$, ie 
$\De\cap\bd B=\bd\De$. 

Suppose \De\ lies in $B_1$. Then it splits $B_1$ into two 3--balls 
one of which, call it $B\shp$, contains $A_0\cup B_0^*$. 
Let $\tau\shp=L\cap B\shp$. 
Then by Lemma 7.4 $(B_0^*,B_0^*\cap L)$ is standardly embedded in 
$(B\shp,B\shp\cap L)$, and by Lemmas 7.3 (4) and 7.5 $(B\shp,B\shp\cap L)$ 
is standardly embedded in $(B_n,B_n\cap L)$ for all $n\geq1$.  

Now suppose that \De\ lies in $\wtM-\inte B_1$. Attaching a 2--handle 
to $B_1$ with core \De\ gives a 3--manifold $Q$ with $\bd Q$ a disjoint 
union of two 2--spheres. One of them bounds a 3--ball $B\shp$ in \wtM\ which 
contains $Q$. $B\shp$ is the union of $B_1$ and a 3--ball $B^+$ such that 
$B\shp\cap B^+$ is a disk. Since by Lemma 7.3 (4) $(B_1,B_1\cap L)$ is 
standardly embedded in 
$(B_{m\shp},B_{m\shp}\cap L)$ for some $m\shp\geq1$ such that 
$B\shp\sbs\inte B_{m\shp}$ we have that $B^+\cap L=\ns$. It follows that 
$(B_0^*, B_0^*\cap L)$ is standardly embedded in $(B\shp,B\shp\cap L)$ which 
by Lemma 7.3 (4) is itself standardly embedded in $(B_n,B_n\cap L)$ for all 
$n\geq m\shp$. 

Note that in both cases we may isotop $B\shp$ slightly so that $\bd B\shp$ 
has at least one fewer intersection curve with $p\n(D)$ than did $\bd B_1$. 
Continuing in this manner we get a 3--ball $B\p$ in \wtM\ such that 
$B\p\cap p\n(D)=\ns$, $(A_0\cup B_0^*)\sbs\inte B\p$, 
$B\p\sbs\inte p\n(C_1\p)$, $(B_0^*,B_0^*\cap L)$ is standardly embedded in 
$(B\p,B\p\cap L)$, and for some $m\p\geq1$ we have that $(B\p,B\p\cap L)$ 
is standardly embedded in $(B_n,B_n\cap L)$ for all $n\geq m\p$. 

We then cut the next 1--handle. We eventually get a 3--ball $B_1^*$ in \wtM\ 
such that $(A_0\cup B_0^*)\sbs\inte B_1^*$, 
$B_1^*\sbs p\n(C_1^*\cap V)\cap\wtV$, $(B_0^*,B_0^*\cap L)$ is standardly 
embedded in $(B_1^*,B_1^*\cap L)$, and for some $m^*\geq 1$ we have that 
$(B_1^*,B_1^*\cap L)$ is standardly embedded in $(B_n,B_n\cap L)$ for 
all $n\geq m^*$. 

We then repeat this argument inductively to construct the exhaustion 
$\{B_n^*\}$. \end{proof}

\section{End reductions at links}

In this section we consider end reductions at non-connected subsets of $M$. 
Let $J$ be a compact but not necessarily connected 3--manifold in $M$. 
As before we define $J$ to be regular if $M-J$ is irreducible and has 
no component with compact closure. End irreducibility and the engulfing  
property rel $J$ in $M$ and end reductions of $M$ at $J$ are defined as 
before. The standard end reduction and uniqueness up to non-ambient 
isotopy work as before. However, $V$ may not be connected, and a component 
of $V\p$ of $V$ may not be an end reduction at $V\p\cap J$. 
As before, each component of $V$ is $\pi_1$--injective in $M$. 

A \textit{regular link} in $M$ is a disjoint union \ga\ of finitely many 
simple closed curves in $M$ such that $M-\ga$ is irreducible. An 
\textit{end reduction of $M$ at \ga } is defined to be an end reduction of 
$M$ at a regular neighborhood of \ga. 

Note that if each component of a link \ga\ is homotopically non-trivial 
in $M$, then \ga\ is a regular link. 

\begin{thm}[Agol] Let $M$ be a connected, orientable, irreducible open 
3--manifold. Let \ga\ be a link in $M$ such that no component of \ga\ is 
homotopically trivial in $M$. Let $V$ be an end reduction of $M$ at 
\ga. Suppose $\pi_1(M)$ is finitely generated 
and decomposable, and that some component of \ga\ is algebraically 
disk busting. Then $V$ is connected and $\pi_1(V)\rightarrow\pi_1(M)$ 
is an isomorphism. \end{thm}

We will prove a slightly more general statement. Let $\ga=\ga_1\cup\cdots
\cup\ga_m$ be a link in $M$ such that no component of \ga\ is homotopically 
trivial in $M$. It is \textit{algebraically disk busting} if  
there is no proper free factorization $A_1*A_2$ of $\pi_1(M)$ 
such that for each $i$ there is a $j$ such that $\ga_i$ is conjugate 
into $A_j$. In general having an algebraically 
disk busting component is sufficient, but not necessary, for \ga\ to be 
algebraically disk busting. If $\pi_1(M)$ is indecomposable, then \ga\ is 
automatically algebraically disk busting.     

\begin{thm} Let $M$ be a connected, orientable, irreducible open 
3--manifold. Let \ga\ be a link in $M$ no component of which is homotopically 
trivial in $M$. Let $V$ be an end reduction 
of $M$ at \ga. Suppose $\pi_1(M)$ is 
finitely generated and \ga\ is algebraically disk busting. Then $V$ is 
connected and $\pi_1(V)\rightarrow\pi_1(M)$ is an isomorphism. \end{thm} 

\begin{proof} First suppose that $\pi_1(M)$ is indecomposable. As in the 
proof of Theorem 4.1 we isotop a compact core $N$ so that it lies in $V^*$. 
Suppose there are distinct components $U$ and $U\p$ of $V^*$ such that 
$U$ is a component of $V$ and $U\p$ contains $N$. Then $\pi_1(U)$ is 
non-trivial. However, as before, a homotopy between a loop in $U$ and a 
loop in $N$ could be cut off on the co-core of a 1--handle that is cut 
in the process of constructing $V^*$. Thus $N$ lies in $V$, and $V$ is 
connected. Hence $\pi_1(V)\rightarrow\pi_1(M)$ is an isomorphism. 

Now suppose that $\pi_1(M)$ is decomposable. As in the proof of Theorem 
5.1 we assume that a compact core $N$ is contained in $\inte C_1$ and 
that the link \ga\ is freely homotopic in $\inte C_1$ to a link $\ga\p$ in 
$\inte N$. We again let $D$ be the co-core of the first 1--handle which is 
cut from $C_1$ and may assume that $N\cap D$ consists of a minimal number 
of essential disks in $N$. 

We may assume that $\ga\p\cap D$ has a minimal number of points. If 
$\ga\p\cap D=\ns$, then some essential disk in $N$ misses $\ga\p$, hence 
the components of $\ga\p$ lie in components of $N$ split along this disk, 
and so the components of \ga\ are conjugate into factors of a proper free 
factorization of $\pi_1(M)$. This contradicts the fact that \ga\ is 
algebraically disk busting. Thus $\ga\p\cap D\neq\ns$. 

Let $S$ be a disjoint union of $m$ copies $S_i^1$ of $S^1$. We have a map 
$f\co S\times[0,1]\rightarrow C_1$ with $f(S\times\{0\})=\ga$ and 
$f(S\times\{1\})=\ga\p$. The argument now proceeds exactly as in 
Theorem 5.1 to contradict either the minimality of $\ga\p\cap D$ or 
that of $N\cap D$.  

Thus $N\cap D=\ns$. We repeat the argument to get that $N\sbs C_1^*$. 
As in the indecomposable case we get that $V$ is connected and $N\sbs V$. 
\end{proof}

\end{document}